\documentclass{amsart}

\usepackage[latin1]{inputenc}
\usepackage{fancyhdr}
\usepackage{graphicx}
\usepackage{newlfont}
\usepackage{amssymb}
\usepackage{amsmath,amscd}
\usepackage{latexsym}
\usepackage{amsthm}
\usepackage[vcentermath]{youngtab}
\usepackage{multirow}
\theoremstyle{plain}                    
\newtheorem{thm}{Theorem}[section]
                          
\newtheorem{claim}[thm]{Claim}

\newtheorem{lem}[thm]{Lemma}
\theoremstyle{definition}
\newtheorem{defn}[thm]{Definition}
\newtheorem{rem}[thm]{Remark}
\newtheorem{exa}[thm]{Example}

\newcommand{\gl}{\operatorname{GL}}
\newcommand{\df}{\displaystyle\frac}
   
\newcommand{\F}{\mathcal{F}}
\newcommand{\End}{\operatorname{End}}
\newcommand{\Ts}{\mathcal{T}^{\mu}}
\newcommand{\T}{\mathcal{T}}
\newcommand{\MM}[2]{M^{#1,#2}(\mathbb{Z}_{\geq 0})}
\newcommand{\card}{\operatorname{card}}
\newcommand{\std}{\operatorname{std}}

\begin{document}
\title[RSK correspondence in the geometry of partial flag varieties]{Robinson-Schensted-Knuth correspondence in the geometry of partial flag varieties}
\author{Daniele Rosso}
\address{\newline
Daniele Rosso \newline 
The University of Chicago \newline 
Department of Mathematics \newline
5734 S. University Ave. Chicago, IL 60615}

\email{d\_rosso@math.uchicago.edu}
\date{\today}
\begin{abstract}
In this paper we generalize to the case of partial flags a result proved both by Spaltenstein and by Steinberg that relates the relative position of two complete flags and the irreducible components of the flag variety in which they lie, using the Robinson-Schensted-Knuth correspondence.
\end{abstract}
\maketitle
\section{Introduction}
\subsection{}
The Robinson-Schensted-Knuth correspondence (RSK for short) is a very classical result. It was first discovered by Robinson (see \cite{R}) as a bijection between permutations of $d$ letters and pairs of standard Young tableaux of the same shape on $d$ boxes, then independently rediscovered by Schensted (see \cite{Sc}). It was eventually generalized by Knuth (see \cite{K}) to the case of two rowed arrays in lexicographic order (or equivalenty matrices with nonnegative integer entries) and pairs of semistandard Young tableaux of the same shape.

This correspondence comes up when considering flag varieties. 
The Bruhat decomposition tells us that the relative position of two complete flags in a $d$-dimensional space, i.e. the orbit of the diagonal action of $\gl_d$ on the variety of pairs of flags in which they lie, is given by an element of the symmetric group $S_d$. 
Now consider the subvariety of flags that are stabilized by a nilpotent transformation $x$. If we let $B$ be a fixed Borel subgroup, we can also see this subvariety as the fixed points of the action of the unipotent element $u=1+x$ on $\gl_d/B$. The irreducible components of this subvariety are parametrized by the standard tableaux on the shape $\lambda$, which is the Jordan type of $x$ (see \cite[II 5.21]{Sp}, \cite{St}). 

Then it is a theorem (see \cite[II 9.8]{Sp}, \cite{St}) that, for two generic flags, their relative position is given by the permutation that we get applying the RSK correspondence to the standard tableaux associated to the irreducible components in which they lie.

\subsection{}
What changes when we look at partial flag varieties instead?

For a nilpotent transformation $x$, we consider the following variety of $n$-step partial flags that are preserved by $x$:
$$\{F:\quad 0=F_0\subset F_1\subset\ldots\subset F_{n-1}\subset F_n=\mathbb{C}^d |x(F_i)\subset F_{i-1} \quad\forall i\}.$$
The irreducible components of this variety can be parametrized by 'semistandard' tableaux (better, by transposed of semistandard tableaux, more on this later) by applying to our specific case some results of Haines about the fibers of convolution morphisms in the affine Grassmanian (see \cite{H}).
This parametrization is also essentially the same that Spaltenstein shows in \cite{Sp1}.

Notice that Shimomura has also worked on partial flag varieties and in \cite{Sh} has given a parametrization of the irreducible components of the variety of partial flags that are invariant under a nilpotent transformation, using Young tableaux, but the variety he considers is different from ours. 

The relative position of two partial flags is given (see for example \cite{blm}) by a matrix with nonnegative integer entries, therefore it seems natural to ask if the theorem generalizes to the case of partial flags. Given two partial flags, is the matrix of relative position the one that corresponds through the more general RSK correspondence to the two semistandard tableaux indexing the irreducible components in which the flags lie? 

As we prove in this paper, the answer is yes, if we modify slightly the usual conventions for the RSK correspondence. We need a variation to account for the fact that the 'semistandard' tableaux mentioned earlier are actually transposed of semistandard tableaux (i.e. the strictness of the inequalities is switched from rows to columns and viceversa).
\subsection*{Acknowledgments}
The author would like to thank Victor Ginzburg for his help and for many ideas and suggestions that led to this paper, and Jonah Blasiak for useful discussions on the subject. In addition, he would also like to thank Joel Kamnitzer for pointing out the result in \cite{H} and Anthony Henderson for pointing out the reference \cite{Sp1}. Finally, he is  grateful to the University of Chicago for support.
 
\section{Flag Varieties and Tableaux}
Let us start by fixing the notation and introducing the objects we will be dealing with. We let $V$ be a $d$-dimensional complex vector space and $\F$ be the variety of complete flags in $V$. We let $x\in \End(V)$ be a nilpotent linear transformation of Jordan type $\lambda=(\lambda_1,\lambda_2,\ldots,\lambda_k)$. Then $\lambda$ is a partition of $d$, which means that  it satisfies $\lambda_1\geq\lambda_2\geq\ldots\lambda_k$, and $|\lambda|=\lambda_1+\lambda_2+\ldots+\lambda_k=d$. 

We consider the subvariety $\F_x\subset\F$ of flags preserved by $x$, that is
$$\F_x:=\{F\in\F|x(F_i)\subset F_{i-1}\}.$$
\begin{defn}\label{deft}
Now let $\T_\lambda$ be the set of standard Young tableaux of shape $\lambda$, we can define a map
$$t:\F_x\to\T_\lambda$$
in the following way: given $F\in\F_x$, consider the Jordan type of the restriction $x|_{F_i}$. This gives us an increasing sequence of Young diagrams each with one box more than the previous one. Filling the new box with the number $i$ at each step, we get a standard tableau.
\end{defn}
Then (see \cite[II 5.21]{Sp},\cite{St}) for a tableau $T\in\T_\lambda$, if we let $\F_{x,T}=t^{-1}(T)\subset\F_x$, we have that the closure $C_T=\overline{\F_{x,T}}$ is an irreducible component of $\F_x$. All the irreducible components are parametrized in this way by the set of standard tableaux of shape $\lambda$. In \cite{Sp}, Spaltenstein actually uses a slightly different parametrization, to see how the two parametrizations are related, see \cite{vL}.
\begin{defn}In this paper, whenever we will refer to a \emph{generic} element in a variety or subvariety, we will mean any element in a suitable open dense subset. \end{defn}
 
We can now state the result ( \cite[II 9.8]{Sp} and \cite[1.1]{St}) that we wish to generalize in this paper.

\begin{thm}\label{stein}Let $\F$ be the variety of complete flags on a vector space $V$, and $x\in\End(V)$ a nilpotent transformation of Jordan type $\lambda$. Let $T,S$ be standard Young tableaux of shape $\lambda$ and $C_T$ and $C_S$ the corresponding irreducible components of $\F_x$. Then for generic flags $F\in C_T$ and $G\in C_S$, the permutation $w(F,G)$ that gives the relative position of the two flags is the same as the permutation $w(T,S)$ given by the RSK correspondence.
\end{thm}

Our goal is to extend this result to varieties of partial flags.

\subsection{Partial Flags and Semistandard Tableaux}\label{pfsst}

Let us fix an integer $n\geq1$ and let $\mu$ be a sequence $\mu=(\mu_1,\ldots,\mu_n)$ of positive integers, such that $|\mu|=\mu_1+\mu_2+\ldots+\mu_n=d$ ($\mu$ is not necessarily a partition because we do not require it to be decreasing). We have the variety of $n$-step flags of type $\mu$ in $V$
$$\F^\mu:=\{F=(0=F_0\subset F_1\subset\ldots\subset F_{n-1}\subset F_n=V)|\dim(F_i/F_{i-1})=\mu_i\}.$$
Then for $x$ as before, we consider the subvariety of partial flags that are preserved by $x$:
$$\F^\mu_x:=\{F\in\F^\mu|x(F_i)\subset F_{i-1}\}.$$
If $F\in\F^\mu_x$, we can associate to $F$ a tableau in an analogous way to definition \ref{deft}, except this time at each step we are adding several boxes, none of which will be in the same row. We start with $1$'s and then at each step we fill in the boxes we just added, with the next integer. The result will be a tableau which is strictly increasing along rows and weakly increasing down columns. For the purpose of this discussion, we will call this kind of tableaux \emph{semistandard}, although by the usual definition this is the transposed of a semistandard tableau.

\begin{defn}Given any tableau $T$ with entries in $\{1,\ldots,n\}$, we say that its \emph{content} is the sequence $\mu=\mu(T)=(\mu_1,\ldots,\mu_n)$ where $\mu_i$ is the number of times the entry $i$ appears in $T$.\end{defn}
 
\begin{defn}\label{eqdeft}
So, if we let $\Ts_\lambda$ be the set of semistandard tableaux of shape $\lambda$ and content $\mu$, we just defined a map
$$t:\F_x^\mu\longrightarrow \Ts_\lambda.$$
\end{defn}

As in the case of standard tableaux, this gives us a bijection between irreducible components of the partial flag variety and the set of semistandard tableaux.
\begin{lem}The irreducible components of $\F_x^\mu$ are the closures $C_T=\overline{\F_{x,T}}$ where $T\in\Ts_\lambda$ and $\F_{x,T}=t^{-1}(T)$.
\end{lem}
For a proof of this fact, see \cite{Sp1} or \cite{H}. Spaltenstein discusses this very briefly, and uses a slightly different convention, as was also mentioned earlier. In his result the indexing set is a subset of the standard tableaux. It can be seen that this subset consists of what we will define later in this paper to be the \emph{standardization} of the semistandard tableaux. 

On the other hand Haines, during the proof of Theorem 3.1 proves a more general result about irreducible components of fibers of convolution morphisms from convolution product of $G(\mathcal{O})$-orbits in the affine Grassmannian. In his result, the combinatorial data are sequences of dominant weights such that the difference of two consecutive weights is in the orbit of the Weyl group acting on a dominant minuscule weight. In our case these correspond to the semistandard tableaux.

\subsection{Relative Position, Words and Arrays} 
Given two flags $F$, $G$ (partial or complete) we define the matrix of nonnegative integers $M(F,G)$ with entries given by:
$$M(F,G)_{ij}=\dim\left(\frac{F_i\cap G_j}{F_i\cap G_{j-1}+F_{i-1}\cap G_j}\right).$$
If $F\in\F^\mu$ and $G\in\F^\nu$, the column sums of this matrix will be $\mu=(\mu_1,\ldots,\mu_n)$ and the row sums will be $\nu=(\nu_1,\ldots,\nu_m)$. 

Then, see \cite[1.1]{blm}, the set $\MM{\mu}{\nu}$ of all such matrices parametrizes the orbits of the  diagonal action of $\gl_d$ on $\F^\mu\times\F^\nu$.
\begin{defn}We call $M(F,G)$ the \emph{relative position} of $F$ and $G$. \end{defn}

In particular, if $F$ and $G$ are both complete flags in $V$, $M(F,G)$ will be a permutation matrix. This data is equivalent to the word $w(F,G)=w(1)\ldots w(d)$  where $w(i)=j$ if $1$ appears in the $(i,j)$-entry of the matrix. 

Similarly, if $F\in\F^\mu$ and $G$ is a complete flag, $M(F,G)$ will have a $1$ in each column and $0$'s everywhere else. Then we can consider the word $w(F,G)=w(1)\ldots w(d)$ defined in the same way. In this case, $w(F,G)$ will not be a permutation word, but a word with content $\mu$.

\begin{defn}\label{arrtom}If $F,G$ are both partial flags, then $M(F,G)$ is just a matrix of nonnegative integers. We can record the same data in a two-rowed array 
$$\omega=\left(\begin{array}{cccc}u_1 & u_2 & \ldots & u_d \\ v_1 & v_2 & \ldots & v_d \end{array}\right)$$
which is defined as follows. 

A pair $\left(\begin{smallmatrix} i \\ j\end{smallmatrix}\right)$ appears in $\omega$ a number of times equal to the $(i,j)$-entry of $M(F,G)$. 

The array $\omega$ is then ordered so that it satisfies the following relation: 
\begin{equation}\label{order}u_1\leq u_2\leq\ldots\leq u_d\qquad\text{ and }\quad v_{k}\geq v_{k+1}\quad\text{ if }u_k=u_{k+1}.\end{equation}
\end{defn}
\begin{exa}\label{notlex}If $M(F,G)$ is the matrix on the left, the corresponding array $\omega$ is given on the right:
$$M(F,G)=\left(\begin{array}{ccc} 1 & 0 & 2 \\ 3 & 1 & 0 \end{array}\right)\quad \omega=\left(\begin{array}{ccccccc} 1 & 1 & 1 & 2 & 2 & 2 & 2 \\ 3 & 3 & 1 & 2 & 1 & 1 & 1 \end{array}\right).$$ \end{exa}
Let $\MM{\mu}{\nu}$ be the set of matrices of non negative integers with row sums $\nu$ and column sums $\mu$. Using the convention just described, we will also identify $\MM{\mu}{\nu}$ with the set of two rowed arrays such that the first row has content $\nu$, the second row has content $\mu$, and they satisfy the order \eqref{order}. 

Depending on what is more convenient at each time, we will use either description of this set.
\begin{rem}\label{lex}Our convention is different from what is used in \cite{F} and \cite{S2}, where the arrays are taken to be in \emph{lexicographic order},
that is with
$$u_1\leq u_2\leq\ldots\leq u_d\qquad\text{ and }\quad v_{k}\leq v_{k+1}\quad\text{ if }u_k=u_{k+1}.$$
With the lexicographic convention, the matrix of example \ref{notlex} would correspond to the array
$$\omega'=\left(\begin{array}{ccccccc} 1 & 1 & 1 & 2 & 2 & 2 & 2 \\ 1 & 3 & 3 & 1 & 1 & 1 & 2 \end{array}\right). $$ 
\end{rem}
\section{Robinson-Schensted-Knuth Correspondence and Standardization} In this section, we will at first review quickly some definitions and properties of the RSK correspondence, following mainly the conventions of \cite[I]{F} and \cite[7.11]{S2}, to which we refer for more details. Then we will see how to adapt the results to the conventions we are using.
\subsection{Review of RSK}
Just for this review, we will call a tableau \emph{semistandard} if it is weakly increasing along rows and strictly increasing down columns. With this convention, the tableaux we defined in section \ref{pfsst} are transposed of semistandard tableaux.
We will also identify matrices with arrays using the lexicographic order, as in Remark \ref{lex}. 

With increasing generality, the RSK correspondence gives a bijection between permutations and pairs of standard tableaux of same shape, or between two-rowed arrays in lexicographic order and pairs of semistandard tableaux of same shape.

Given a permutation word $w=w(1)\ldots w(d)$ or more generally a two rowed array $\omega=\left(\begin{smallmatrix}u_1 & u_2 & \ldots & u_d \\ v_1 & v_2 & \ldots & v_d \end{smallmatrix}\right)$, the algorithm is given by inserting the entries of the word (or of the second row of the array) by row bumping in the first tableau. At the same time we need to record in the second tableau which box has been added at each step (in the more general case of the array, the added box at the $k$-th step will be recorded with $u_k$ as opposed to $k$). The convention for row bumping is that a new entry $z$ bumps the left-most entry in the row which is \emph{strictly larger} than $z$.

If $T$, $S$ are semistandard tableaux and $\omega$ is an array in lexicographic order, we will denote the correspondence by
$$ M(T,S)=\omega; \quad\text{ or }\quad(T,S)\stackrel{RSK}{\longleftrightarrow}\omega.$$

As can be seen in \cite[7.11]{S2}, given a semistandard tableau $T$ we can consider its \emph{standardization} $\tilde{T}$. It is a standard tableau of the same shape as $T$, constructed in this way: the $\mu_1$ boxes that contain $1$ in $T$ will be replaced by the numbers $1,2,\ldots,\mu_1$ increasingly from left to right. Then the boxes that originally contained $2$'s will be replaced by $\mu_1+1,\ldots,\mu_1+\mu_2$ always increasingly from left to right, and so on.
\begin{exa}\label{ex1}
$$T=\young(112,2,3)\qquad\tilde{T}=\young(124,3,5).$$
\end{exa}
In a similar way, given an array in lexicographic order $\omega=\left(\begin{smallmatrix}u_1 & u_2 & \ldots & u_d \\ v_1 & v_2 & \ldots & v_d \end{smallmatrix}\right)$ we can define the standardization $\tilde{\omega}$. It is given by replacing $u_i$ with $i$ in the first row, while in the second row we replace the $1$'s with $1,2,\ldots,\mu_1$ increasing from left to right, then the $2$'s and so on. Notice that the standardization of an array is a permutation.
\begin{exa}\label{ex2}
$$\omega=\left(\begin{array}{ccccc} 1 & 2 & 2 & 3 & 3 \\ 3 & 1 & 2 & 1 & 2 \end{array}\right)\qquad \tilde{\omega}=\left(\begin{array}{ccccc} 1 & 2 & 3 & 4 & 5 \\ 5 & 1 & 3 & 2& 4\end{array}\right).$$
\end{exa}
Then the important result is that standardization allows us to always reduce the RSK correspondence to the special case of permutations and standard tableaux, because standardization and RSK commute.
\begin{lem}\label{lem1}The following diagram commutes:
$$\begin{CD}
\Ts_\lambda\times\T_\lambda^{\nu} @>\text{RSK}>> M^{\mu,\nu}(\mathbb{Z}_{\geq 0}) \\
@VV\std\times\std V                             @VV\std V \\
\T_\lambda\times\T_\lambda @>\text{RSK}>>S_d
\end{CD}$$
\end{lem}
In the diagram, $\T_\lambda$, $\Ts_\lambda$, $\T_\lambda^{\nu}$ are respectively the set of standard tableaux and the sets of semistandard tableaux with content $\mu$ and $\nu$, all of shape $\lambda$. Also, $S_d$ is the set of permutations of $d$ elements and $M^{\mu,\nu}(\mathbb{Z}_{\geq 0})$ is is the set of two rowed arrays in lexicographic order with the first row having content $\nu$ and the second row having content $\mu$ (identified with matrices). Finally, $\std$ is the standardization map.

The lemma is proved in \cite[7.11.6]{S2}, but let us illustrate this with an example. 
\begin{exa}Let $T$, $\omega$ be as in examples \ref{ex1} and \ref{ex2} and let 
$$S=\young(123,2,3)\quad\text{ then we have }\quad \std(S)=\tilde{S}=\young(135,2,4)$$
then $(T,S)\stackrel{RSK}{\longleftrightarrow}\omega$
and indeed $(\tilde{T},\tilde{S})\stackrel{RSK}{\longleftrightarrow}\tilde{\omega}$.
\end{exa}
\subsection{Variation on RSK}\label{varrsk} In this paper we will need a slight variation on the RSK correspondence. This will agree with RSK on permutations, but will give different results in the case of general two rowed arrays. In particular, it will associate to a two rowed array satisfying \eqref{order}, a pair of tableaux that are strictly increasing along rows and weakly increasing down columns. This is what we called \emph{semistandard} in section \ref{pfsst} and we will keep using this terminology from now on. In the rest of this paper, we will also set the convention of identifying matrices and arrays using Definition \ref{arrtom}.

The variation of the correspondence is defined modifying the row bumping algorith to the following: a new entry $z$ will bump the left-most entry in the row which is \emph{greater or equal} to $z$. The recording tableau will be then constructed in the usual way. 

This difference is clearly irrelevant in the case of standard tableaux, but note that our new choice of row bumping will in fact produce tableaux that are strictly increasing along rows and weakly increasing down columns. This is similar to the dual RSK defined in \cite[7.14]{S2}, which however is only defined for matrices of $0$'s and $1$'s.

Since we will only use this variation on the correspondence, from now on we will call this one RSK and we will use the same notation as before, there should be no confusion.

\begin{lem}This procedure gives a bijection between matrices of non-negative integers and pairs of semistandard (strictly increasing along rows and weakly increasing down columns) tableaux of same shape.\end{lem} 
\begin{proof}This is completely analogous to the usual proofs of the RSK correspondence (see \cite{F},\cite{S2}).

If the array corresponding to the matrix is $\omega=\left(\begin{array}{ccc}u_1 & \ldots & u_d \\ v_1 & \ldots & v_d \end{array}\right)$ and by the correspondence it gives us the pair of tableaux $(P,Q)$, then it is clear that the insertion tableau $P$ will be semistandard. To check that the recording tableau $Q$ is also semistandard, it is enough to show that if $u_i=u_{i+1}$, then $u_{i+1}$ will end up in a row of $Q$ that is strictly below the row of $u_i$. 

Since $\omega$ satisfies \eqref{order}, if $u_i=u_{i+1}$, then $v_i\geq v_{i+1}$. This means that if $v_i$ bumps an element $y_i$ from the first row, then the element $y_{i+1}$ bumped by $v_{i+1}$ from the first row must be in the same box where $y_i$ was or in a box to the left of it. In turn, this implies that $y_i\geq y_{i+1}$ and we can iterate this argument for the following rows. Now, the bumping route $R_{i}$ of $v_{i}$ must stop before the bumping route $R_{i+1}$ of $v_{i+1}$, which will then continue at least one row below that of $R_i$, which shows what we want.  

The fact that the correspondence is a bijection just follows from the fact that we can do the reverse row bumping algorithm by taking at each step the box that in the recording tableau contains the biggest number. In case of equal elements, we will take the one that is in the lowest row.
\end{proof}
\begin{rem}
Basically in this version of RSK we are considering equal entries in a tableau to be 'bigger' if they are in a lower row and, while inserting, sequences of equal numbers are considered decreasing sequences.\end{rem}

This leads us to a new definition of \emph{standardization} that will give us an analogous result to lemma \ref{lem1}. 
Given a semistandard tableau $T$, we define its standardization $\tilde{T}$ by replacing the $1$'s with $1,2,\ldots,\mu_1$ \emph{starting from the top row and going down}, and then the same for $2$'s and so on. For an array $\omega$ ordered as in \eqref{order}, we define $\tilde{\omega}$ by replacing the first row with $1,2,\ldots,d$ and on the second row we replace the $1$'s by $1,2,\ldots,\mu_1$ \emph{decreasingly} from left to right and same for the rest, always decreasing from left to right.
\begin{exa}\label{exsst}
$$T=\young(12,12,3)\qquad \tilde{T}=\young(13,24,5)$$
$$\omega=\left(\begin{array}{ccccc} 1 & 2 & 2 & 3 & 3 \\ 1 & 3 & 1 & 2 & 2 \end{array}\right)\qquad
\tilde{\omega}=\left(\begin{array}{ccccc} 1 & 2 & 3 & 4 & 5 \\ 2 & 5 & 1 & 4 & 3 \end{array}\right) $$
\end{exa}
Now, with our new conventions for semistandard tableaux, order of arrays, RSK and standardization and the same notation of lemma \ref{lem1} we have that
\begin{lem}\label{lem2}Standardization and RSK commute, as in the following diagram:
$$\begin{CD}
\Ts_\lambda\times\T^{\nu}_\lambda @>\text{RSK}>> M^{\mu,\nu}(\mathbb{Z}_{\geq 0}) \\
@VV\std\times\std V                             @VV\std V \\
\T_\lambda\times\T_\lambda @>\text{RSK}>>S_d
\end{CD}$$
\end{lem}
The proof, mutatis mutandis, is the same as the proof of lemma \ref{lem1} in \cite[7.11.6]{S2}. It is just the observation that the standardization we choose for the arrays is exactly the one that makes the insertion procedure work the way we want, turning sequences of equal numbers into decreasing sequences.
\begin{exa}Let $T$, $\omega$ as in example \ref{exsst} and let
$$S=\young(12,23,3)\quad\text{ then we have }\quad \tilde{S}=\young(12,34,5)$$
then $(T,S)\stackrel{RSK}{\longleftrightarrow}\omega$ and $(\tilde{T},\tilde{S})\stackrel{RSK}{\longleftrightarrow}\tilde{\omega}$.
\end{exa}
\begin{rem}\label{rmk}It is very easy to see that if we fix the contents $\mu$ and $\nu$, two different arrays $\omega_1\neq\omega_2\in M^{\mu,\nu}(\mathbb{Z}_\geq 0)$ when standardized will give two different permutations $\tilde{\omega}_1\neq\tilde{\omega}_2$. That is we have an injective map
$$\std:\MM{\mu}{\nu}\to S_d.$$
We therefore have an inverse
$$\std^{-1}:\std(\MM{\mu}{\nu})\to \MM{\mu}{\nu}$$
which is easily described as follows:
$$\left(\begin{array}{cccccc} 1 & 2 & \ldots & \nu_1 & \nu_1 +1 & \ldots \\ v_1 & v_2 & \ldots & v_{\nu_1} & v_{\nu_1+1} & \ldots \end{array}\right)
\mapsto \left(\begin{array}{cccccc} 1 & 1 & \ldots & 1 & 2 & \ldots \\ v_1' & v_2' & \ldots & v_{\nu_1}' & v_{\nu_1+1}' & \ldots\end{array}\right)
 $$
the first row is just replaced by $\nu_1$ $1$'s, followed by $\nu_2$ $2$'s and so on, while we have 
$$v_k'=j\quad\text{ if }\quad v_k\in\{\mu_1+\ldots+\mu_{j-1}+1,\ldots,\mu_1+\ldots+\mu_j\}.$$ 
\end{rem}
\section{RSK and Partial Flag Varieties} In this section we will use all the conventions of section \ref{varrsk} and the notations of section \ref{pfsst}.

We can now state and prove the main result, which generalizes Theorem \ref{stein}. The strategy for the proof is to use standardization and Lemma \ref{lem2} to reduce the problem to the case of complete flags.
\begin{thm}Let $x\in\End(V)$ be a nilpotent transformation of Jordan type $\lambda$, $T\in\Ts_\lambda$, $S\in\T^{\nu}_\lambda$ be semistandard tableaux, and let $C_T$ and $C_S$ be respectively the irreducible components of $\F^\mu_x$ and $\F^\nu_x$ corresponding to the tableaux $T$ and $S$. 

Then, for generic $F\in C_T$ and $G\in C_S$, we have that the relative position matrix $M(F,G)$ is the same as the matrix $M(T,S)$ given by the RSK correspondence. \end{thm}
\begin{proof}
For a fixed $\mu=(\mu_1,\ldots,\mu_n)$ with $|\mu|=\mu_1+\ldots+\mu_n=d$, consider the map 
$$p_\mu:\F\to\F^\mu$$ 
that forgets some of the spaces, that is
$$(0=F_0,F_1,F_2,\ldots,F_{n-1},F_n=V)\mapsto (0=F_0,F_{\mu_1},F_{\mu_1+\mu_2},\ldots,F_{\mu_1+\ldots+\mu_{n-1}},F_n=V).$$

Clearly, if $F$ is any partial flag in $\F^\mu_x$ and $\tilde{F}\in p_\mu^{-1}(F)$, then $\tilde{F}\in\F_x$ because for all $j$ there is some some $i$ such that
$$\tilde{F}_{\mu_1+\ldots+\mu_i}\subset\tilde{F}_{j-1}\subset\tilde{F}_j\subset\tilde{F}_{\mu_1+\ldots+\mu_{i+1}}$$
and 
$$x(\tilde{F}_j)\subset x(\tilde{F}_{\mu_1+\ldots+\mu_{i+1}})=x(F_{i+1})\subset F_i=\tilde{F}_{\mu_1+\ldots+\mu_i}\subset\tilde{F}_{j-1}.$$

Now, let $t:\F_x^\mu\to\Ts_\lambda$ be the map that associates a semistandard tableau to a partial flag, as in Definition \ref{eqdeft}. 

We fix a semistandard tableau $T$ and we let $\F_{x,T}:=t^{-1}(T)$, then $\F_{x,T}$ is a constructible dense subset of $C_T$.

Let $\tilde{T}$ be the standardization of $T$ and let $\F_{x,\tilde{T}} =t^{-1}(\tilde{T})\subset\F_x$ be the dense subset of $C_{\tilde{T}}$. The set $C_{\tilde{T}}$ is the irreducible component of the complete flag variety associated to the standard tableau $\tilde{T}$.

It is clear that if $\tilde{F}\in \F_{x,\tilde{T}}$, then we have $F=p_\mu(\tilde{F})\in \F_{x,T}$ because
$$x|_{F_i}=x|_{\tilde{F}_{\mu_1+\ldots+\mu_i}}$$
also, the map
$$p_\mu:\F_{x,\tilde{T}}\to \F_{x,T}$$
is surjective. This is because we can always find appropriate subspaces to complete a partial flag $F$ to a flag $\tilde{F}$ such that the restriction of $x$ to those subspaces has the Jordan type we want.

What we have said so far applies in the same way if we fix a semistandard tableau $S$ of content $\nu$ and we consider the sets $\F_{x,S}\subset\F^\nu_x$ and $\F_{x,\tilde{S}}\subset\F_x$.

Now, let us fix two semistandard tableaux $T$ and $S$ as in the statement of the theorem, and consider their standardizations $\tilde{T}$ and $\tilde{S}$. For generic complete flags $\tilde{F}\in C_{\tilde{T}}$ and $\tilde{G}\in C_{\tilde{S}}$, Theorem \ref{stein} tells us that $M(\tilde{F},\tilde{G})=M(\tilde{T},\tilde{S})$. We let then $X_{\tilde{T}}\subset C_{\tilde{T}}$, $X_{\tilde{S}}\subset C_{\tilde{S}}$ be the open dense subsets such that this is true.

Then $X_{\tilde{T}}\cap \F_{x,\tilde{T}}$ is constructible dense in $C_{\tilde{T}}$. Hence it contains an open dense subset and the image of
$$p_\mu:X_{\tilde{T}}\cap \F_{x,\tilde{T}}\to \F_{x,T}$$
is constructible dense in $\F_{x,T}$, therefore it is also dense in $C_T$. In the same way, $p_\nu(X_{\tilde{S}}\cap \F_{x,\tilde{S}})$ is constructible dense in $\F_{x,S}$.
\begin{claim}If $F\in p_\mu(X_{\tilde{T}}\cap \F_{x,\tilde{T}})$ and $G\in p_\nu(X_{\tilde{S}}\cap \F_{x,\tilde{S}})$ then $M(F,G)=M(T,S)$. \end{claim}
Let $\tilde{F}\in p_\mu^{-1}(F)$ and $\tilde{G}\in p_\nu^{-1}(G)$, then by Lemma \ref{lem2} we have that 
$$\tilde{\omega}=\std(M(T,S))=M(\tilde{T},\tilde{S}).$$

Now let $\omega'=M(F,G)$. By the definition of relative position of flags, the array $\tilde{\omega}=M(\tilde{T},\tilde{S})=M(\tilde{F},\tilde{G})$ is such that for all $i,j$
$$ \card \left\{\begin{tabular}{r|l}
\multirow{2}{*}{$\left(\begin{array}{c} \tilde{u} \\ \tilde{v} \end{array}\right)\in \tilde{\omega}$}
& $\tilde{u}\in\{\nu_1+\ldots+\nu_{j-1}+1,\ldots,\nu_1+\ldots\nu_j\}$, \\
& $\tilde{v}\in\{\mu_1+\ldots+\mu_{i-1}+1,\ldots,\mu_1+\ldots+\mu_i\}$
\end{tabular}\right\} $$

\begin{align*} &= \dim\left(\df{\tilde{F}_{\mu_1+\ldots+\mu_i}\cap\tilde{G}_{\nu_1+\ldots+\nu_j}}{(\tilde{F}_{\mu_1+\ldots+\mu_{i-1}}\cap \tilde{G}_{\nu_1+\ldots+\nu_j}) + (\tilde{F}_{\mu_1+\ldots+\mu_{i}}\cap \tilde{G}_{\nu_1+\ldots+\nu_{j-1}})}\right) \\
 &= \dim\left(\df{F_i\cap G_j}{F_i\cap G_{j-1}+F_{i-1}\cap G_j}\right) \\
 &= \card\left\{\left(\begin{array}{c} u \\ v \end{array}\right)\in\omega' \left| \left(\begin{array}{c} u \\ v \end{array}\right)\right.= \left(\begin{array}{c} i \\ j \end{array}\right)\right\}
 \end{align*}
Therefore, by Remark \ref{rmk}, $\omega'=\std^{-1}(\tilde{\omega})$. It follows that $\std(\omega')=\tilde{\omega}$, that is 
$$\std(M(F,G))=\std(M(T,S)).$$ 
Again by Remark \ref{rmk}, this implies that $M(F,G)=M(T,S)$. This concludes the proof of the claim. 

Since $p_\mu(X_{\tilde{T}}\cap \F_{x,\tilde{T}})$ and $p_\nu(X_{\tilde{S}}\cap \F_{x,\tilde{S}})$ are constructible dense in $C_T$ and $C_S$ respectively, they each contain an open dense subset of the respective irreducible component, which proves the theorem.
\end{proof}

\end{document}